\documentclass[11pt,reqno]{amsart}
\usepackage{amssymb,amscd}
\newtheorem{thm}{Theorem}[section]
\newtheorem{prop}[thm]{Proposition}
\newtheorem{lemma}[thm]{Lemma}

\def\U{\mathrm{U}}
\def\bC{\mathbb{C}}
\def\cM{\mathcal{M}}
\def\bN{\mathbb{N}}
\def\eps{\varepsilon}
\def\tr{\mathrm{tr}}
\def\proj{\mathrm{proj}}
\def\cN{\mathcal{N}}
\def\1{\mathbf{1}}
\def\cA{\mathcal{A}}
\def\cB{\mathcal{B}}
\def\ffi{\varphi}

\def\Ric{\mathrm{Ric}}
\def\u{\mathfrak{u}}
\def\h{\mathfrak{h}}
\def\g{\mathfrak{g}}
\def\<{\langle}
\def\>{\rangle}
\def\Re{\mathrm{Re}}
\def\Tr{\mathrm{Tr}}
\def\bZ{\mathbb{Z}}

\begin{document}

\title[A log-Sobolev type inequality]{A log-Sobolev type inequality for free
entropy \\ of two projections}
\author[F. Hiai]{Fumio Hiai$\,^{1,2}$}
\address{Graduate School of Information Sciences,
Tohoku University, Aoba-ku, Sendai 980-8579, Japan}
\author[Y. Ueda]{Yoshimichi Ueda$\,^{1,3}$}
\address{Graduate School of Mathematics,
Kyushu University, Fukuoka 810-8560, Japan}

\thanks{$^1\,$Supported in part by Japan Society for the Promotion of Science,
Japan-Hungary Joint Project.}
\thanks{$^2\,$Supported in part by Grant-in-Aid for Scientific Research
(B)17340043.}
\thanks{$^3\,$Supported in part by Grant-in-Aid for Young Scientists
(B)17740096.}

\maketitle

\begin{abstract}
We prove an inequality between the free entropy and the mutual free Fisher
information for two projections, regarded as a free analog of the logarithmic
Sobolev inequality. The proof is based on the random matrix approximation
procedure via the Grassmannian random matrix model of two projections.
\end{abstract}

\section*{Introduction}

Among the most important notions in classical information theory is the mutual
information formally expressed as
\begin{equation}\label{F-0.1}
I(X,Y)=H(X)+H(Y)-H(X,Y)
\end{equation}
for two random variables $X,Y$ in terms of their Shannon-Gibbs entropies
$H(\cdot)$. Motivated by the above expression, in part VI \cite{V6} of his series
of papers, Voiculescu introduced the notions of the mutual free Fisher information
$\ffi^*$ and of the mutual free information $i^*$ via the liberation theory in
free probability. These quantities are defined for subalgebras (rather than random
variables) of a tracial $W^*$-probability space while the microstates free entropy
$\chi$ as well as the non-microstates $\chi^*$ is for random variables. In the
last section of \cite{V6}, Voiculescu explained that the formula like \eqref{F-0.1}
\begin{align}\label{F-0.2}
&\chi(X_1,\dots,X_n)+\chi(Y_1,\dots,Y_m) \nonumber\\
&\qquad=i^*(W^*(X_1,\dots,X_n),W^*(Y_1,\dots,Y_m))
+\chi(X_1,\dots,X_n,Y_1,\dots,Y_m)
\end{align}
cannot be true in general, and he suggested the necessity of generalizing the free
entropy to more general objects beyond self-adjoint variables and proposed, 
for instance, how to define the free entropy for projections. Note here that 
the free entropy $\chi$ for self-adjoint variables has always $-\infty$ for
projections.

On the other hand, a large deviation principle recently obtained in
\cite{HP1} is  related to a pair $(P(N),Q(N))$ of random projection matrices
having independent and unitarily invariant distribution provided $\lim_N
\mathrm{rank}(P(N))/N$ and $\lim_N\mathrm{rank}(Q(N))/N$ exist as constants.
Indeed, the random pair $(P(N),Q(N))$ induces a random tracial state $\tau_N$ on
$C^*(\bZ_2\star\bZ_2)$, the universal $C^*$-algebra generated by two
projections, and  the large deviation principle is concerned with the empirical 
measure of $\tau_N$. An important fact there is that the rate function
$\mathcal{I}(\tau)$ of tracial  states $\tau$ is equal to the minus of the free
entropy $\chi_\proj(p,q)$ of two  projection generators $(p,q)$ in the GNS
representation with respect to $\tau$.

The main aim of this paper is to prove the inequality
\begin{equation}\label{F-0.3}
-\chi_\proj(p,q)\le\ffi^*(p:q)
\end{equation}
for two projections $(p,q)$ under mild assumptions, where $\ffi^*(p:q)$ is the
mutual free Fisher information of subalgebras $\bC p+\bC(\1-p)$ and
$\bC q+\bC(\1-q)$. The proof is based on a random matrix approximation procedure
derived from the large deviation principle mentioned above. In fact, the inequality
\eqref{F-0.3} arises as a scaling limit of the classical logarithmic Sobolev
inequality due to Bakry and Emery \cite{BE} applied  to a Grassmannian random
matrix pair modeled on the pair $(p,q)$. Thus we may consider \eqref{F-0.3} as a
kind of free probabilistic logarithmic Sobolev inequality. Such free analogs have
been previously obtained in \cite{Bi} for single self-adjoint variables and in
\cite{HPU} for single unitary variables. After then a remarkable approach to such
free analogs is given in \cite{L} for single self-adjoint variables.  

The paper is organized as follows. First in \S1, we briefly recall the definitions
of the free entropy $\chi_\proj$ for projections and of the mutual free Fisher
information $\ffi^*$. For convenience of reference, the explicit forms of
$\chi_\proj(p,q)$ in \cite{HP1} and of $\ffi^*(p:q)$ in \cite{V6} for two
projections are mentioned. \S2 is devoted to the proof of \eqref{F-0.3} based on
the classical logarithmic Sobolev inequality in \cite{BE} via the Grassmannian 
random matrix approximation. Here we need the Ricci curvature tensor of the 
Grassmannian manifold to verify Bakry and Emery's $\Gamma_2$-criterion. \S3
contains supplementary remarks. We note that $-\chi_\proj(p,q)$ appears as a 
scaling limit of the classical mutual information on the Grassmannian manifold;
it seems natural because
$-\chi_\proj(p,q)=\chi_\proj(p)+\chi_\proj(q)-\chi_\proj(p,q)$ (due to
$\chi_\proj(p) = \chi_\proj(q) = 0$) has the form like \eqref{F-0.1}. From the
viewpoint in \cite{V6}, this form also suggests that $-\chi_\proj(p,q)$ should
coincide with the mutual free information $i^*(p,q)$. In fact, we make a heuristic
computation to indicate that $i^*(p,q)=-\chi_\proj(p,q)$, a very particular case
of \eqref{F-0.2}, since $\chi_\proj(p)=0$ for any single projection $p$. Finally,
the equality \eqref{F-0.3} is slightly generalized into the relative version
including a certain potential term.

\medskip\noindent
{\it Acknowledgment.} The authors thank Prof.~Masaki Izumi for fruitful discussions.    

\section{Preliminaries}
\setcounter{equation}{0}

\subsection{Free entropy for projections}
For $N\in\bN$ let $\U(N)$ be the unitary group of order $N$. For $k\in\{0,1,\dots,N\}$
let $G(N,k)$ denote the set of all $N\times N$ orthogonal projection matrices of
rank $k$, that is, $G(N,k)$ is identified with the so-called Grassmannian manifold
consisting of $k$-dimensional subspaces in $\bC^N$. Let $P_N(k)$ be the diagonal
matrix with the first $k$ diagonals $1$ and the others $0$. Each $P\in G(N,k)$ is
diagonalized as
\begin{equation}\label{F-1.1}
P=UP_N(k)U^*,
\end{equation}
where  $U\in\U(N)$ is determined modulo $\U(k)\oplus\U(N-k)$. Hence $G(N,k)$ is
identified with the homogeneous space $\U(N)/(\U(k)\oplus\U(N-k))$, and the unitarily
invariant probability measure on $G(N,k)$ corresponds to the measure on
$\U(N)/(\U(k)\oplus\U(N-k))$ induced from the Haar probability measure
$\gamma_{\U(N)}$ on $\U(N)$. We denote by $\gamma_{G(N,k)}$ this unitarily invariant
measure on $G(N,k)$. Let $\zeta_{N,k}:\U(N)\to G(N,k)$ be the (surjective continuous)
map defined by the equation \eqref{F-1.1}, i.e., $\zeta_{N,k}(U):=UP_N(k)U^*$. Then
the measure $\gamma_{G(N,k)}$ is more explicitly written as
\begin{equation}\label{F-1.2}
\gamma_{G(N,k)}=\gamma_{\U(N)}\circ\zeta_{N,k}^{-1}.
\end{equation}

Let $(p_1,\dots,p_n)$ be an $n$-tuple of projections in a tracial $W^*$-probability
space $(\cM,\tau)$ with $\alpha_i:=\tau(p_i)$, $1\le i\le n$. Following Voiculescu's
proposal in \cite[14.2]{V6} we define the {\it free entropy}
$\chi_\proj(p_1,\dots,p_n)$ of $(p_1,\dots,p_n)$ as follows: Choose
$k(N,i)\in\{0,1,\dots,N\}$ for each $N\in\bN$ and $1\le i\le n$ in such a way that
$k(N,i)/N\to\alpha_i$ as $N\to\infty$ for $1\le i\le n$. For each $m\in\bN$ and
$\eps>0$ set
\begin{eqnarray*}
&&\Gamma_\proj(p_1,\dots,p_n;k(N,1),\dots,k(N,n);N,m,\eps) \\
&&\quad:=\biggl\{(P_1,\dots,P_n)\in\prod_{i=1}^nG(N,k(N,i)):
\big|\tr_N(P_{i_1}\cdots P_{i_r})-\tau(p_{i_1}\cdots p_{i_r})\big|<\eps
\nonumber\\
&&\hskip6.5cm\mbox{for all $1\le i_1,\dots,i_r\le n$, $1\le r\le m$}\biggr\},
\end{eqnarray*}
where $\tr_N$ stands for the normalized trace on the $N\times N$ matrices. We then
define
\begin{eqnarray}\label{F-1.3}
&&\chi_\proj(p_1,\dots,p_n) 
:=\inf_{m\in\bN,\,\eps>0}\limsup_{N\to\infty} \nonumber\\
&&\quad{1\over N^2}\log
\Biggl(\bigotimes_{i=1}^n\gamma_{G(N,k(N,i))}\Biggr)
\bigl(\Gamma(p_1,\dots,p_n;k(N,1),\dots,k(N,n);N,m,\eps)\bigr).
\end{eqnarray}
It is easy to see that the above definition of $\chi_\proj(p_1,\dots,p_n)$ is
independent of the choices of $k(N,i)$ with $k(N,i)/N\to\alpha_i$ for $1\le i\le n$.
The free entropy $\chi_\proj$ for projections has properties similar to those for
self-adjoint variables developed in \cite{V1,V2,V4} and for unitary variables in
\cite[Chapter 6]{HP}, which we will discuss elsewhere \cite{HU}. It is obvious that
$\chi_\proj(p)=0$ for any single projection $p$.

In this paper we are concerned with the free entropy of two projections. Let $(p,q)$
be a pair of projections in $(\cM,\tau)$ with $\alpha:=\tau(p)$ and $\beta:=\tau(q)$.
Set
$$
E_{11}:=p\wedge q,\quad E_{10}:=p\wedge q^\perp,\quad
E_{01}:=p^\perp\wedge q,\quad E_{00}:=p^\perp\wedge q^\perp,
$$
$$
E:=\1-(E_{00}+E_{01}+E_{10}+E_{11})
$$
and $\alpha_{ij}:=\tau(E_{ij})$ for $i,j=0,1$. Then $E$ and $E_{ij}$ are in the
center of $\cN:=\{p,q\}''$ and $(E\cN E,\tau|_{E\cN E})$ is isomorphic to
$L^\infty((0,1),\nu;M_2(\bC))$, the $L^\infty$-algebra of $M_2(\bC)$-valued
functions, where $\nu$ is a measure on $(0,1)$ with
$\nu((0,1))=1-\sum_{i,j=0}^1\alpha_{ij}$. Here $EpE$ and $EqE$ correspond to
$$
t\in(0,1)\mapsto\bmatrix1&0\\0&0\endbmatrix
\ \ \mbox{and}\ \ \bmatrix t&\sqrt{t(1-t)}\\\sqrt{t(1-t)}&1-t\endbmatrix,
$$
respectively, and $\tau|_{E\cN E}$ is represented as
\begin{equation}\label{F-1.4}
\tau(a)=\int_0^1\tr_2(a(t))\,d\nu(t)
\end{equation}
for $a\in E\cN E$ corresponding to $a(\cdot)\in L^\infty((0,1),\nu;M_2(\bC))$.
In this way, the mixed moments of $(p,q)$ with respect to $\tau$ are determined
by the data $(\nu,\{\alpha_{ij}\}_{i,j=0}^1)$. Although $\nu$ is not necessarily
a probability measure, we define the free entropy $\Sigma(\nu)$ by
$$
\Sigma(\nu):=\int_0^1\int_0^1\log|x-y|\,d\nu(x)\,d\nu(y)
$$
in the same fashion as in \cite{V1}. Furthermore, we set
\begin{equation}\label{F-1.5}
\rho:=\min\{\alpha,\beta,1-\alpha,1-\beta\}
=\frac{1}{2}\Biggl(1-\sum_{i,j=0}^1 \alpha_{ij}\Biggr),
\end{equation}
\begin{equation}\label{F-1.6}
C:=\rho^2B\biggl({|\alpha-\beta|\over\rho},{|\alpha+\beta-1|\over\rho}\biggr)
\end{equation}
(meant zero if $\rho=0$), where
\begin{align*}
B(s,t)&:={(1+s)^2\over2}\log(1+s)-{s^2\over2}\log s
+{(1+t)^2\over2}\log(1+t)-{t^2\over2}\log t \\
&\qquad-{(2+s+t)^2\over2}\log(2+s+t)+{(1+s+t)^2\over2}\log(1+s+t)
\end{align*}
for $s,t\geq0$. With these definitions, the following expression was obtained in
\cite{HP1} as a consequence of the large deviation principle for an independent
pair of random projection matrices.

\begin{prop}\label{P-1.1}{\rm (\cite[Theorem 3.2, Proposition 3.3]{HP1})}\quad
If $\alpha_{00}\alpha_{11}=\alpha_{01}\alpha_{10}=0$, then
\begin{align*}
\chi_\proj(p,q)&={1\over4}\Sigma(\nu)+{\alpha_{01}
+\alpha_{10}\over2}\int_0^1\log x\,d\nu(x) \\
&\qquad\qquad+{\alpha_{00}+\alpha_{11}\over2}\int_0^1\log(1-x)\,d\nu(x)-C,
\end{align*}
and otherwise $\chi_\proj(p,q)=-\infty$.
\end{prop}

It is known \cite{HP1} that $\limsup$ in definition \eqref{F-1.3} can be replaced by
$\lim$ in the case of two projections. Furthermore, it was shown there that
$\chi_\proj(p,q)=0$ if and only if $p$ and $q$ are free. (Note that this fact still
remains valid even for general $n$ tuples of projections, whose proof will be given
in \cite{HU}.) Note that the condition
$\alpha_{00}\alpha_{11}=\alpha_{01}\alpha_{10}=0$ is equivalent to
\begin{equation}\label{F-1.7}
\begin{cases}
\alpha_{11}=\max\{\alpha+\beta-1,0\}, \\
\alpha_{00}=\max\{1-\alpha-\beta,0\}, \\
\alpha_{10}=\max\{\alpha-\beta,0\}, \\
\alpha_{01}=\max\{\beta-\alpha,0\};
\end{cases}
\end{equation}
in this case, $\alpha_{01}+\alpha_{10}=|\alpha-\beta|$ and
$\alpha_{00}+\alpha_{11}=|\alpha+\beta-1|$.

\subsection{Mutual free Fisher information}
Let $\cA$ ($\ni\1$) and $\cB$ ($\ni\1$) be two $*$-subalgebras in $(\cM,\tau)$,
which are assumed to be algebraically free. Let $\cA\vee\cB$ and $W^*(\cA\cup\cB)$
denote the subalgebra and the von Neumann subalgebra, respectively, generated by
$\cA\cup\cB$. Let $\delta_{\cA:\cB}$ be the derivation from $\cA\vee\cB$ into the
$\cA\vee\cB$-bimodule $(\cA\vee\cB)\otimes(\cA\vee\cB)$ uniquely determined by
$$
\begin{cases}
\delta_{\cA:\cB}(a)=a\otimes\1-\1\otimes a & \text{for $a\in\cA$}, \\
\delta_{\cA:\cB}(b)=0 & \text{for $b\in\cB$}.
\end{cases}
$$
If there is an element $\xi\in L^1(W^*(\cA\cup\cB))$ such that
$$
\tau(\xi x)=(\tau\otimes\tau)(\delta_{\cA:\cB}(x)),\qquad x\in\cA\vee\cB,
$$
then $\xi$ is called the {\it liberation gradient} of $(\cA,\cB)$ and denoted by
$j(\cA:\cB)$. Voiculescu \cite{V6} introduced the {\it mutual free Fisher information}
of $\cA$ relative to $\cB$ by
$$
\ffi^*(\cA:\cB):=\|j(\cA:\cB)\|_2^2
$$
($\|\cdot\|_2$ stands for the $L^2$-norm with respect to $\tau$) if $j(\cA:\cB)$
exists in $L^2(W^*(\cA\cup\cB))$; otherwise $\ffi^*(\cA:\cB):=+\infty$. See \cite{V6}
for more about the mutual free Fisher information.

Let $(p,q)$ be a pair of projections in $(\cM,\tau)$ and set $\cA:=\bC p+\bC(\1-p)$,
$\cB:=\bC q+\bC(\1-q)$. Then the liberation gradient $j(\cA:\cB)$ and the mutual free
Fisher information $\ffi^*(\cA:\cB)$ were computed in \cite{V6}. Here, recall that
the Hilbert transform of a function $f$ with $f(x)/(1+|x|) \in L^1(\mathbb{R},dx)$
is defined to be 
$$
(Hf)(x):=\lim_{\eps\searrow0}(H_{\varepsilon}f)(x) \quad \text{with}
\quad (H_{\varepsilon}f)(x) := \int_{|x-t|>\eps}{f(t)\over x-t}\,dt
$$
(whenever the limit exists almost everywhere). The following is a slightly
improved version of \cite[Proposition 12.7]{V6}:  

\begin{prop}\label{P-1.2} With the same notations as in \S\S1.1 assume that
$\alpha_{00}\alpha_{11}=\alpha_{01}\alpha_{10}=0$, $\nu$ has the density
$f:=d\nu/dx \in L^3((0,1),x(1-x)dx)$ and moreover
\begin{equation}\label{integrability}
\int_0^1 \biggl(\frac{\alpha_{01}+\alpha_{10}}{x}
+ \frac{\alpha_{11}+\alpha_{00}}{1-x}\biggr) f(x)\,dx<+\infty.
\end{equation}
Define $X:=pqp+(\1-p)(\1-q)(\1-p)$ and
$$
\phi(x):=(Hf)(x)+{\alpha_{01}+\alpha_{10}\over x}
-{\alpha_{00}+\alpha_{11}\over1-x}\quad\mbox{for $0<x<1$}.
$$
Then
$$
j(\cA:\cB)=[q,p]\phi(EXE)\in L^2(\cM,\tau)
$$
and hence
$$
\ffi^*(\cA:\cB)=\int_0^1\phi(x)^2f(x)x(1-x)\,dx<+\infty.
$$
\end{prop}

The assumption \eqref{integrability} can be reduced when $\alpha=\beta$ or
$\alpha+\beta = 1$ thanks to \eqref{F-1.7}. In fact, \eqref{integrability} is nothing
when $\alpha=\beta=1/2$; it means $\int_0^1 x^{-1}f(x)\,dx < +\infty$ when
$\alpha+\beta = 1$ but $\alpha\neq\beta$. All the assumptions of Proposition
\ref{P-1.2} are satisfied in particular when $p$ and $q$ are free (see
\cite[Example 3.6.7]{VDN}). Note \cite[Propositions 5.17 and 9.3.c]{V6} that
$j(\cA:\cB)=0$ (or equivalently $\ffi^*(\cA:\cB)=0$) if and only if $p$ and $q$ are
free.

In \cite[\S12]{V6} the support of $\nu$ was assumed to be an infinite set to
guarantee that $\cA$ and $\cB$ are algebraically free. As long as $\nu\ne0$, that is
automatically satisfied from the assumption of $\nu$ having the density. In the case
where $\nu=0$ so that $\rho=0$ by \eqref{F-1.5}, it follows that $p\in\{0,\1\}$ or
$q\in\{0,\1\}$; hence Proposition \ref{P-1.2} trivially holds.

\medskip\noindent
{\it Proof of Proposition 1.2.}\enspace
We first remark a weighted norm version of so-called M.~Riesz's theorem for Hilbert
transform. Since $\int_0^1 (x(1-x))^{-1/2}\,dx <+\infty$, the celebrated weighted
norm inequality for Hilbert transform \cite[Theorem 8]{Hunt at el} shows that there
is a constant $C_w>0$ depending only on the weight function
$w(x):=\1_{[0,1]}(x)x(1-x)$ such that for every function $g$ (whose $Hg$ can be
defined)    
\begin{equation}\label{zzz} 
\Vert Hg \Vert_{w,3} \leq C_w\Vert g\Vert_{w,3}
\end{equation}  
with the weighted norm
$$
\Vert g \Vert_{w,p} := \biggl(\int_0^1 |g(x)|^p x(1-x)\,dx\biggr)^{1/p}
\quad\mbox{for $1\le p<\infty$},
$$
and moreover $\Vert H_{\varepsilon} g - Hg\Vert_{w,3} \rightarrow 0$ as
$\varepsilon\searrow0$ whenever $\Vert g \Vert_{w,3} < +\infty$. In what follows we
use the same symbols as in \cite[\S12]{V1} with small exception; $p,q$, $\cA, \cB$
and $x,x_1,x_2$ are used instead of $P,Q$, $A, B$ and $t,t_1,t_2$, respectively. By
the facts mentioned above one easily has
\allowdisplaybreaks{ 
\begin{align*} 
&\int_0^1\int_0^1 \biggl(\frac{x_1^{n+1}-x_2^{n+1}}{x_1-x_2}
- \frac{x_1^n-x_2^n}{x_1-x_2}\biggr)\,d\nu(x_1)\,d\nu(x_2) \\
&\quad=
-\lim_{\varepsilon\searrow0}\iint_{|x_1-x_2|>\eps}\biggl(x_1^{n-1} x_1(1-x_1) \frac{f(x_2)}{x_1-x_2}f(x_1) \\
&\hskip4cm+ x_2^{n-1} x_2(1-x_2)\frac{f(x_1)}{x_2-x_1}f(x_2)\biggr)\,dx_1\,dx_2 \\
&\quad= 
-2\lim_{\varepsilon\searrow0}\int_0^1 x^{n-1} (H_{\varepsilon}f)(x) f(x) x(1-x)\,dx \\
&\quad=
-2\int_0^1 x^{n-1} (Hf)(x) f(x) x(1-x)\,dx.       
\end{align*} 
}Hence, the assertion of \cite[Lemma 12.6]{V6} can be changed to 
\begin{align} 
((\tau\otimes\tau)\circ\delta_{\cB:\cA}((pq)^n) &= 
-\frac{1}{2}\int_0^1 x^{n-1} (Hf)(x) f(x) x(1-x)\,dx \notag \\
&\quad+ 
\frac{1-\alpha}{2}\int_0^1 x^{n-1}(x-1)\,d\nu(x)
+ \frac{\alpha_{00}+\alpha_{11}}{2} \int_0^1 x^{n-1}\,d\nu(x) \label{xxx}
\end{align} 
under the assumptions of Proposition \ref{P-1.2}. The rest of the proof goes along
the same line as \cite[Proposition 12.7]{V6} with replacing \cite[Lemma 12.6]{V6}
by \eqref{xxx}. \qed

\section{An inequality}
\setcounter{equation}{0}

Our aim of this section is to obtain the following inequality between the free
entropy $\chi_\proj$ and the mutual free Fisher information $\ffi^*$ for a pair
$(p,q)$ of projections in a $W^*$-probability space $(\cM,\tau)$. For simplicity we
hereafter write $\ffi^*(p:q)$ for the mutual free Fisher information
$\ffi^*(\bC p+\bC(\1-p):\bC q+\bC(\1-q))$ (see \S\S1.2).

\begin{thm}\label{T-2.1}
With the same assumptions as stated in Proposition \ref{P-1.2},
$$
-\chi_\proj(p,q)\le\ffi^*(p:q).
$$
\end{thm}

The main idea of the proof is a random matrix approximation procedure based on
the large deviation shown in \cite{HP1}. In fact, we apply Bakry and Emery's
logarithmic Sobolev inequality in \cite{BE} to random projection matrix pairs (or
probability measures on the product of two Grassmannian manifolds) and pass to
the scaling limit as the matrix size goes to $\infty$. Thus our inequality can
be regarded as a kind of free probability counterpart of the logarithmic Sobolev
inequality. A further discussion on this aspect will be given in the next section.

So-called Bakry and Emery's $\Gamma_2$-criterion is crucial in their logarithmic
Sobolev inequality in the Riemannian manifold setting and the Ricci curvature
tensor is one of  the important ingredients of the criterion. We thus need to
compute the Ricci curvature tensor $\Ric(G(N,k))$ of $G(N,k)$, $1\le k\le N-1$,
as described below. Let $\u(N)$ be the Lie algebra of $\U(N)$ and regard
$\h(N,k):=\u(k)\oplus\u(N-k)$ as a Lie subalgebra of $\u(N)$. 
The tangent space $T_PG(N,k)$ at each $P\in G(N,k)$ can be identified
with $\g(N,k):=\h(N,k)^\perp$, the orthocomplement of $\h(N,k)$ in $\u(N)$ with
respect to the Riemannian metric $\<X,Y\>:=\Re\,\Tr_N(XY^*)$, where $\Tr_N$ is 
the usual trace on $N\times N$ matrices. Choose the following complete orthonormal
system of $\g(N,k)$:
\begin{equation}\label{F-2.1}
E_{ij}:={1\over\sqrt2}(e_{ij}-e_{ji}),\quad
F_{ij}:={\sqrt{-1}\over\sqrt2}(e_{ij}+e_{ji})
\end{equation}
with $1\le i\le k$, $k+1\le j\le N$. According to well-known facts on compact
matrix groups and O'Neill's formula (see \cite[Proposition 3.17, Theorem 3.61]{GHL}
for example), the Ricci curvature tensor of $G(N,k)$  with respect to the
above-mentioned Riemannian metric is computed as follows:
\begin{eqnarray*}
&&\Ric(G(N,k))_P(X,X) \\
&&\qquad=\sum_{1\le i\le k,\,k+1\le j\le N}
\bigl(\|[X,E_{ij}]\|_{HS}^2+\|[X,F_{ij}]\|_{HS}^2\bigr),
\qquad X\in\g(N,k).
\end{eqnarray*}
A simple direct computation shows that the above right-hand side is $N\|X\|_{HS}^2$
so that
\begin{equation}\label{F-2.2}
\Ric(G(N,k))=NI_{2k(N-k)}.
\end{equation}

\bigskip\noindent
{\it Proof of Theorem 2.1.}\enspace
Let $\alpha$, $\beta$ and $(\nu,\{\alpha_{ij}\}_{i,j=0}^1)$ be as in \S\S1.1 for
the given pair $(p,q)$ of projections. Since the inequality trivially holds if
$\nu=0$, assume $\nu\ne0$ and let $\nu_1:=\nu(1)^{-1}\nu$, the normalization 
of $\nu$. In addition to the assumptions of Proposition \ref{P-1.2} we first assume
the following (A) and (B):
\begin{itemize}
\item[(A)] $\nu$ is supported in $[\delta,1-\delta]$ for some $\delta>0$ and it
has the continuous density $d\nu/dx$.
\item[(B)] The function
$$
Q_{\nu_1}(x):=2\int_0^1\log|x-y|\,d\nu_1(y)
$$
is a well-defined $C^1$-function on $[0,1]$.
\end{itemize}
Choose $C^1$-functions $h_0(x)$ and $h_1(x)$ on $[0,1]$ such that
$$
h_0(x)\begin{cases}
=\log x & \text{($\delta\le x\le1$)}, \\
\ge\log x & \text{($0\le x\le\delta$)},
\end{cases}
\qquad
h_1(x)\begin{cases}
=\log(1-x) & \text{($0\le x\le1-\delta$)}, \\
\ge\log(1-x) & \text{($1-\delta\le x\le1$)}.
\end{cases}
$$
For each $N\in\bN$ choose $k(N),l(N)\in\{1,\dots,N-1\}$ such that
$k(N)/N\to\alpha$ and $l(N)/N\to\beta$ as $N\to\infty$, and set
\begin{align*}
n_0(N)&:=N-\min\{k(N),l(N)\}, \\
n_1(N)&:=\max\{k(N)+l(N)-N,0\}, \\
n(N)&:=\min\{k(N),l(N),N-k(N),N-l(N)\}=N-n_0(N)-n_1(N).
\end{align*}
For each $n\in\bN$ letting
\begin{align*}
\psi_N(x):={n(N)\over N}Q_{\nu_1}(x)&+{|k(N)-l(N)|\over N}\,h_0(x) \\
&+{|k(N)+l(N)-N|\over N}\,h_1(x),\qquad0\le x\le1,
\end{align*}
we define a probability measure (regarded as a pair of $N\times N$ random projection
matrices) $\lambda_N^{\psi_N}$ on $G(N,k(N))\times G(N,l(N))$ by
\begin{equation}\label{F-2.3}
d\lambda_N^{\psi_N}(P,Q):={1\over Z_N^{\psi_N}}
\exp\bigl(-N\Tr_N(\psi_N(PQP))\bigr)
\,d\bigl(\gamma_{G(N,k(N))}\otimes\gamma_{G(N,l(N))}\bigr)(P,Q)
\end{equation}
with the normalization constant $Z_N^{\psi_N}$ as well as the reference measure
$\lambda_N^0:=\gamma_{G(N,k(N))}\otimes\gamma_{G(N,l(N))}$. When
$(P,Q)\in G(N,k(N))\times G(N,l(N))$ is distributed under
$\lambda_N^0$, it is known \cite[(2.1)]{HP1} that the eigenvalues of $PQP$ are
\begin{equation}\label{F-2.4}
\underbrace{0,\dots,0}_{n_0(N)\ {\rm times}},
\underbrace{1,\dots,1}_{n_1(N)\ {\rm times}}, \, x_1,\dots,x_{n(N)}
\end{equation}
and the joint distribution of $(x_1,\dots,x_{n(N)})$ is
\begin{align}\label{F-2.5}
d\tilde\lambda_N^0(x_1,\dots,x_{n(N)})
&:={1\over\widetilde Z_N^0}\prod_{i=1}^{n(N)}
x_i^{|k(N)-l(N)|}(1-x_i)^{|k(N)+l(N)-N|} \nonumber\\
&\qquad\quad\times\prod_{1\le i<j\le n(N)}
(x_i-x_j)^2\prod_{i=1}^{n(N)}\1_{[0,1]}(x_i)\,dx_i
\end{align}
with the normalization constant $\widetilde Z_N^0$. Hence it turns out that when
$(P,Q)\in G(N,k(N))\times G(N,l(N))$ is distributed under $\lambda_N^{\psi_N}$,
the eigenvalues of $PQP$ are listed as in \eqref{F-2.4} but the joint
distribution of $(x_1,\dots,x_{n(N)})$ is changed to
\begin{eqnarray}\label{F-2.6}
&&d\tilde\lambda_N^{\psi_N}(x_1,\dots,x_{n(N)}) \nonumber\\
&&\quad:={1\over\widetilde Z_N^{\psi_N}}
\exp\Biggl(-\sum_{i=1}^{n(N)}\bigl\{n(N)Q_{\nu_1}(x_i)
+|k(N)-l(N)|(h_0(x_i)-\log x_i) \nonumber\\
&&\hskip4cm+|k(N)+l(N)-N|(h_1(x_i)-\log(1-x_i))\bigr\}\Biggr) \nonumber\\
&&\hskip3cm\times\prod_{1\le i<j\le n(N)}(x_i-x_j)^2
\prod_{i=1}^{n(N)}\1_{[0,1]}(x_i)\,dx_i
\end{eqnarray}
with another normalization constant $\widetilde Z_N^{\psi_N}$.

Similarly to \cite[Proposition 2.1]{HP1} and \cite[\S5.5]{HP} we have
\begin{itemize}
\item[(a)] The limit
$C':=\lim_{N\to\infty}{1\over N^2}\log\widetilde Z_N^{\psi_N}$ exists as well as
$C=\lim_{N\to\infty}{1\over N^2}\log\widetilde Z_N^0$ (see \eqref{F-1.6}).
\item[(b)] When $(x_1,\dots,x_{n(N)})$ is distributed under
$\tilde\lambda_N^{\psi_N}$, the empirical measure
$$
{\delta_{x_1}+\dots+\delta_{x_{n(N)}}\over n(N)}
$$
satisfies the large deviation principle in the scale $1/N^2$ with the rate
function
$$
I(\mu):=-\rho^2\Sigma(\mu)+\rho^2\int_0^1F(x)\,d\mu(x)+C'
\quad\mbox{for $\mu\in\cM([0,1])$},
$$
where $\cM([0,1])$ is the set of probability measures on $[0,1]$, $\rho$ is given
in \eqref{F-1.5} and
$$
\quad F(x):=Q_{\nu_1}(x)+{|\alpha-\beta|\over\rho}(h_0(x)-\log x)
+{|\alpha+\beta-1|\over\rho}(h_1(x)-\log(1-x))
$$
for $0\le x\le1$.
\item[(c)] $\nu_1$ is a unique minimizer of $I$ with $I(\nu_1)=0$.
\end{itemize}
The last assertion follows from \cite[I.1.3 and I.3.1]{ST} because by the construction
of $h_0$ and $h_1$ we get
$$
Q_{\nu_1}(x)\begin{cases}
=F(x) & \text{if $x\in[\delta,1-\delta]\supset{\rm supp}\,\nu_1$}, \\
\le F(x) & \text{for $x\in[0,1]$}.
\end{cases}
$$
Furthermore, the above large deviation yields:
\begin{itemize}
\item[(d)] The mean eigenvalue distribution
$$
\hat\lambda_N^{\psi_N}:=\int_{[0,1]^{n(N)}}
{\delta_{x_1}+\dots+\delta_{x_{n(N)}}\over n(N)}
\,d\tilde\lambda_N^{\psi_N}(x_1,\dots,x_{n(N)})
$$
weakly converges to $\nu_1$ as $N\to\infty$.
\end{itemize}

Since the Riemannian manifold $G(N,k(N))\times G(N,l(N))$ has the volume measure
$\lambda_N^0$ and its Ricci curvature tensor is $NI_{2k(N-k)+2l(N-l)}$ by
\eqref{F-2.2}, the classical logarithmic Sobolev inequality due to Bakry and
Emery \cite{BE} implies that
\begin{equation}\label{F-2.7}
S(\lambda_N^{\psi_N},\lambda_N^0)
\le{1\over2N}\int_{G(N,k(N))\times G(N,l(N))}
\bigg\|\nabla\log{d\lambda_N^{\psi_N}\over d\lambda_N^0}\bigg\|_{HS}^2
\,d\lambda_N^{\psi_N},
\end{equation}
where the left-hand side is the relative entropy of $\lambda_N^{\psi_N}$ with
respect to $\lambda_N^0$ and the gradient 
$\nabla\log(d\lambda_N^{\psi_N}/d\lambda_N^0)(P,Q)$ 
is considered in $\mathfrak{g}(N,k(N))\oplus\mathfrak{g}(N,l(N))$ via the natural
identification 
$T_{(P,Q)} G(N,k(N))\times G(N,l(N)) = \mathfrak{g}(N,k(N))\oplus
\mathfrak{g}(N,l(N))$. By \eqref{F-2.5} and \eqref{F-2.6} notice that
$$
{d\lambda_N^{\psi_N}\over d\lambda_N^0}(P,Q)
={1\over Z_N^{\psi_N}}\exp\bigl(-N\Tr_N(\psi_N(PQP))\bigr)
={\widetilde Z_N^0\over\widetilde Z_N^{\psi_N}}
\exp\Biggl(-N\sum_{i=1}^{n(N)}\psi_N(x_i)\Biggr)
$$
for $(P,Q)\in G(N,k(N))\times G(N,l(N))$ and for the eigenvalues
$(x_1,\dots,x_{n(N)})$ of $PQP$ except $n_0(N)$ zeros and $n_1(N)$ ones (see
\eqref{F-2.4}). Hence we get
\begin{align*}
S(\lambda_N^{\psi_N},\lambda_N^0)
&=\int_{G(N,k(N))\times G(N,l(N))}\log
{d\lambda_N^{\psi_N}\over d\lambda_N^0}(P,Q)\,d\lambda_N^{\psi_N}(P,Q) \\
&=\log\widetilde Z_N^0-\log\widetilde Z_N^{\psi_N}
+\int_{[0,1]^{n(N)}}\Biggl(-N\sum_{i=1}^{n(N)}\psi_N(x_i)\Biggr)
\,d\tilde\lambda_N^{\psi_N}(x_1,\dots,x_{n(N)}) \\
&=\log\widetilde Z_N^0-\log\widetilde Z_N^{\psi_N}
-Nn(N)\int_0^1\psi_N(x)\,d\hat\lambda_N^{\psi_N}(x).
\end{align*}
Since $\psi_N(x)$ converges to
$$
\rho Q_{\nu_1}(x)+|\alpha-\beta|h_0(x)+|\alpha+\beta-1|h_1(x)
$$
uniformly on $[0,1]$, it follows from (a) and (d) above that
\begin{eqnarray*}
&&\lim_{N\to\infty}{1\over N^2}S(\lambda_N^{\psi_N},\lambda_N^0) \\
&&\quad=C-C'-\rho\int_0^1\bigl(\rho Q_{\nu_1}(x)+|\alpha-\beta|h_0(x)
+|\alpha+\beta-1|h_1(x)\bigr)\,d\nu_1(x).
\end{eqnarray*}
Since (c) gives
$$
-C'=-\rho^2\Sigma(\nu_1)+\rho^2\int_0^1F(x)\,d\nu_1(x) \\
=-\rho^2\Sigma(\nu_1)+\rho^2\int_0^1Q_{\nu_1}(x)\,d\nu_1(x),
$$
we have
\begin{eqnarray}\label{F-2.8}
&&\lim_{N\to\infty}{1\over N^2}S(\lambda_N^{\psi_N},\lambda_N^0) \nonumber\\
&&\quad=C-\rho^2\Sigma(\nu_1)-\rho\int_0^1\bigl(|\alpha-\beta|h_0(x)
+|\alpha+\beta-1|h_1(x)\bigr)\,d\nu_1(x) \nonumber\\
&&\quad=-\chi_\proj(p,q)
\end{eqnarray}
thanks to Proposition \ref{P-1.1} and \eqref{F-1.7} together with
$\nu_1=(2\rho)^{-1}\nu$ (see \cite[(3.4)]{HP1}).

On the other hand, since
$$
\nabla\log{d\lambda_N^{\psi_N}\over d\lambda_N^0}(P,Q)
=-N\nabla\bigl(\Tr_N(\psi_N(PQP))\bigr),
$$
one can compute
$$
\bigg\|\nabla\log{d\lambda_N^{\psi_N}\over d\lambda_N^0}(P,Q)\bigg\|_{HS}^2
=4N^2\Tr_N\bigl((\psi_N'(PQP))^2PQP(I-PQP)\bigr),
$$
whose short proof will be given as Lemma \ref{L-2.2} below for completeness.
Therefore,
\begin{eqnarray*}
&&\int_{G(N,k(N))\times G(N,l(N))}
\bigg\|\nabla\log{d\lambda_N^{\psi_N}\over d\lambda_N^0}(P,Q)\bigg\|_{HS}^2
\,d\lambda_N^{\psi_N}(P,Q) \\
&&\quad=4N^2\int_{[0,1]^{n(N)}}\sum_{i=1}^{n(N)}
(\psi_N'(x_i))^2x_i(1-x_i)\,d\tilde\lambda_N^{\psi_N}(x_1,\dots,x_{n(N)}) \\
&&\quad=4N^2n(N)\int_0^1(\psi_N'(x))^2x(1-x)\,d\hat\lambda_N^{\psi_N}(x) \\
&&\quad=4n(N)\int_0^1\bigl(n(N)Q_{\nu_1}'(x)+|k(N)-l(N)|h_0'(x) \\
&&\hskip3cm+|k(N)+l(N)-1|h_1'(x)\bigr)^2x(1-x)\,d\hat\lambda_N^{\psi_N}(x),
\end{eqnarray*}
and thus by (d) we have
\begin{eqnarray}\label{F-2.9}
&&\lim_{N\to\infty}{1\over2N^3}\int_{G(N,k(N))\times G(N,l(N))}
\bigg\|\nabla\log{d\lambda_N^{\psi_N}\over d\lambda_N^0}\bigg\|_{HS}^2
\,d\lambda_N^{\psi_N} \nonumber\\
&&\quad=2\rho\int_0^1\bigl(\rho Q_{\nu_1}'(x)+|\alpha-\beta|h_0'(x)
+|\alpha+\beta-1|h_1'(x)\bigr)^2x(1-x)\,d\nu_1(x) \nonumber\\
&&\quad=\int_0^1\biggl(\rho Q_{\nu_1}'(x)+{|\alpha-\beta|\over x}
-{|\alpha+\beta-1|\over1-x}\biggr)^2x(1-x)\,d\nu(x) \nonumber\\
&&\quad=\ffi^*(p:q)
\end{eqnarray}
thanks to Proposition \ref{P-1.2}, since $\nu_1 = (2\rho)^{-1}\nu$ so that
$\rho Q_{\nu_1}'(x)=(Hf)(x)$ for $f:=d\nu/dx$. Combining \eqref{F-2.7}--\eqref{F-2.9}
yields the desired inequality under the additional assumptions (A) and (B).

Next, let us remove (A) and (B). First, suppose that the assumption (A) is still
satisfied but (B) is not. For each $\eps>0$ choose a non-negative $C^\infty$-function
$\psi_\eps$ supported in $[-\eps,\eps]$ with $\int\psi_\eps(x)\,dx=1$. Let
$f_\eps:=f*\psi_\eps$ for $f:=d\nu/dx$ and define $d\nu_\eps(x):=f_\eps(x)\,dx$;
then $\nu_\eps$ is a measure supported in a closed proper subinterval of $(0,1)$
with $\nu_\eps((0,1))=1-\sum_{i,j=0}^1\alpha_{ij}$ whenever $\eps$ is small enough.
Let $(p_\eps,q_\eps)$ be a pair of projections in some $(\cM,\tau)$ corresponding to
the representing data $(\nu_\eps,\{\alpha_{ij}\}_{i,j=0}^1)$. (Such a pair can be
constructed via the GNS representation of the universal $C^*$-algebra
$C^*(\bZ_2\star\bZ_2)$ with respect to the tracial state corresponding to
$(\nu_\eps,\{\alpha_{ij}\}_{i,j=0}^1)$; see \cite[\S3]{HP1} and also \S\S3.3.) Since
(A) and (B) are satisfied for $\nu_\eps$, we get
$$
-\chi_\proj(p_\eps,q_\eps)\le\ffi^*(p_\eps:q_\eps).
$$
Since $\|f_\eps-f\|_{w,3}\to0$ and $\|Hf_\eps-Hf\|_{w,3}\to0$ as $\eps\searrow0$ (see the proof of Proposition 1.2 for the weighted norm $\| \cdot \|_{w,3}$),
the H\"older inequality together with \eqref{zzz} implies that
\begin{equation}\label{conv}
\int_0^1 ((Hf_{\eps})(x))^2 f_{\eps}(x) x(1-x)\,dx 
\longrightarrow \int_0^1 ((Hf)(x))^2 f(x) x(1-x)\,dx
\end{equation}
and hence
$$
\lim_{\eps\searrow0}\ffi^*(p_\eps:q_\eps)=\ffi^*(p:q).
$$
Since $\Sigma(\mu)$ for $\mu\in\cM((0,1))$ is weakly upper semicontinuous (see
\cite[5.3.2]{HP}), we also have
$$
-\chi_\proj(p,q)\le\liminf_{\eps\searrow0}
\bigl(-\chi_\proj(p_\eps,q_\eps)\bigr)
$$
so that $-\chi_\proj(p,q)\le\ffi^*(p:q)$.

Finally, suppose only the assumptions stated in Proposition \ref{P-1.2}. For
$\delta>0$ set
$$
d\nu_\delta(s):={1-\sum_{i,j=0}^1\alpha_{ij}\over\nu([\delta,1-\delta])}
\1_{[\delta,1-\delta]}(x)\,d\nu(x)
$$
and let $(p_\delta,q_\delta)$ be a pair of projections corresponding to
$(\nu_\delta,\{\alpha_{ij}\}_{i,j=0}^1)$. Let us denote the density of $\nu_{\delta}$
by $f_\delta$; then it is immediate to see that $\|f_\delta-f\|_{w,3}\to0$. To show
that $\ffi^*(p_\delta:q_\delta) \rightarrow \ffi^*(p:q)$ as $\delta\searrow0$,
it suffices to prove the following convergences as $\delta\searrow0$:    
\begin{align}
\int_0^1 ((Hf_{\delta})(x))^2 f_{\delta}(x) x(1-x)\,dx 
&\longrightarrow \int_0^1 ((Hf)(x))^2 f(x) x(1-x)\,dx, \label{2.1-1}\\
\int_0^1 (Hf_{\delta})(x)x^{-1}f_{\delta}(x)\,x(1-x)\,dx 
&\longrightarrow \int_0^1 (Hf)(x) x^{-1}f(x)\,x(1-x)\,dx, \label{2.1-2}\\
\int_0^1 (Hf_{\delta})(x) (1-x)^{-1}f_{\delta}(x)\,x(1-x)\,dx 
&\longrightarrow \int_0^1 (Hf)(x)(1-x)^{-1}f(x)\,x(1-x)\,dx, \label{2.1-3}\\
\int_0^1 x^{-2}f_{\delta}(x)\,x(1-x)\,dx 
&\longrightarrow \int_0^1 x^{-2}f(x)\,x(1-x)\,dx, \label{2.1-4}\\
\int_0^1 (1-x)^{-2}f_{\delta}(x)\,x(1-x)\,dx 
&\longrightarrow \int_0^1 (1-x)^{-2}f(x)\,x(1-x)\,dx. \label{2.1-5}   
\end{align} 
Remark here that \eqref{2.1-2} and \eqref{2.1-4} are unnecessary when
$\alpha_{01}+\alpha_{10} = |\alpha-\beta|=0$, and so are \eqref{2.1-3} and
\eqref{2.1-5} when $\alpha_{11}+\alpha_{00}=|\alpha+\beta-1|=0$. 
The convergence \eqref{2.1-1} follows as \eqref{conv} above. Also, \eqref{2.1-4} and
\eqref{2.1-5} immediately follow from the hypothesis \eqref{integrability}. Since
\eqref{2.1-2} and \eqref{2.1-3} are similarly shown, let us prove only the former
here. Thus, we should assume $\alpha\neq\beta$, and \eqref{integrability} means
$\int_0^1 x^{-1}f(x)\,dx < +\infty$. By the H\"older inequality together with
\eqref{zzz} one can estimate  
\begin{align*} 
\|(Hf_{\delta})x^{-1}f_{\delta}-(Hf)x^{-1}f\|_{w,1} 
&\leq 
\|H(f_{\delta}-f)\|_{w,3}\cdot \|x^{-1}f_{\delta}^{1/2}\|_{w,2}\cdot\|f_{\delta}^{1/2}\|_{w,6} \\ 
&\quad+ \|Hf\|_{w,3}\cdot\|x^{-1}|f_{\delta}-f|^{1/2}\|_{w,2}\cdot\|\,|f_{\delta}-f|^{1/2}\|_{w,6} \\
&\leq 
C_w\|f_{\delta}-f\|_{w,3}\cdot \|x^{-1}f_{\delta}^{1/2}\|_{w,2}\cdot\|f_{\delta}\|_{w,3}^{1/2} 
\\ 
&\quad+ C_w\|f\|_{w,3}\cdot\|x^{-1}|f_{\delta}-f|^{1/2}\|_{w,2}\cdot\|f_{\delta}-f\|_{w,3}^{1/2}.   
\end{align*} 
Note that 
$$
\|x^{-1}f_{\delta}^{1/2}\|_{w,2}^2 
\le\int_0^1 x^{-1}f_{\delta}(x)\,dx
\longrightarrow \int_0^1 x^{-1}f(x)\,dx,
$$
$$
\|x^{-1}|f_{\delta}-f|^{1/2}\|_{w,2}^2 
\le \int_0^1 x^{-1}|f_{\delta}(x)-f(x)|\,dx
\longrightarrow 0
$$
as $\delta\searrow0$, where $\int_0^1 x^{-1}f(x)\,dx < +\infty$ is essential. These
apparently imply \eqref{2.1-2} thanks to $f \in L^3((0,1),x(1-x)dx)$ and
$\|f_{\delta}-f\|_{w,3} \longrightarrow0$.  

Moreover, since $-\log x<x^{-1}$ near $0$ and $-\log(1-x)<(1-x)^{-1}$ near $1$, the hypothesis \eqref{integrability} implies that 
$$
\int_0^1(-\log x)f_\delta(x)\,dx
\longrightarrow\int_0^1(-\log x)f(x)\,dx<+\infty,
$$
$$
\int_0^1(-\log(1-x))f_\delta(x)\,dx
\longrightarrow\int_0^1(-\log(1-x))f(x)\,dx<+\infty
$$
as $\delta\searrow0$ (whenever those are needed) so that   
$$
-\chi_\proj(p,q)\le\liminf_{\delta\searrow0}
\bigl(-\chi_\proj(p_\delta,q_\delta)\bigr). 
$$
Hence the proof is completed. \qed

\begin{lemma}\label{L-2.2}
Let $\psi$ be a $C^1$-function on $[0,1]$ and define $\Psi(P,Q):=\Tr_N(\psi(PQP))$
for $(P,Q)\in G(N,k)\times G(N,l)$. Then
$$
\|\nabla\Psi(P,Q)\|_{HS}^2=4\Tr_N\bigl((\psi'(PQP))^2PQP(I-PQP)\bigr)
$$
holds for every $(P,Q)\in G(N,k)\times G(N,l)$.
\end{lemma}

\proof
Write $(X_r)_{r=1}^{2k(N-k)}$ for the orthonormal basis of ${\frak g}(N,k)$ given in
\eqref{F-2.1} and also $(Y_s)_{s=1}^{2l(N-l)}$ for that of ${\frak g}(N,l)$. For each
$(P,Q)=(UP_N(k)U^*,VP_N(l)V^*)$ in $G(N,k)\times G(N,l)$, a local normal coordinate
at $(P,Q)$ is given by the mapping
\begin{align}\label{F-2.10}
&(X,Y)\in{\frak g}(N,k)\oplus{\frak g}(N,l) \notag \\ 
&\qquad\mapsto(Ue^XP_N(k)e^{-X}U^*,Ve^YP_N(l)e^{-Y}V^*)\in G(N,k)\times G(N,l).
\end{align}
By a direct computation using this coordinate, one can compute
\begin{align*}
\nabla\Psi(P,Q)&=\sum_r\bigl\<U^*QPf'(PQP)PU-U^*Pf'(PQP)PQU,X_r\bigr\>X_r \\
&\quad+\sum_s\bigl\<V^*Pf'(PQP)PQV-V^*QPf'(PQP)PV,Y_s\bigr\>Y_s
\end{align*}
so that
\begin{align*}
\|\nabla\Psi(P,Q)\|_{HS}^2
&=2\|P_k(N)U^*Pf'(PQP)PQU(I-P_k(N))\|_{HS}^2 \\
&\quad+2\|P_l(N)V^*QPf'(PQP)PV(I-P_l(N))\|_{HS}^2 \\
&=4\Tr_N\bigl((f'(PQP))^2PQP(I-PQP)\bigr).
\end{align*}
\qed

\section{Supplementary remarks}
\setcounter{equation}{0}

\subsection{Classical vs free probabilistic mutual information}
The classical mutual information of two random variables, say $X, Y$, is usually
formulated to be 
$$
I(X,Y) := S(\mu_{(X,Y)},\mu_X\otimes\mu_Y)
= \int_{\mathcal{X}\times\mathcal{X}}
\log\frac{d\mu_{(X,Y)}}{d(\mu_X\otimes \mu_Y)}(x,y)\,d\mu_{(X,Y)}(x,y),
$$
where $\mu_X,\mu_Y$ are the distribution measures of $X,Y$ on the phase space
$\mathcal{X}$ and $\mu_{(X,Y)}$ the joint distribution of $(X,Y)$ on
$\mathcal{X}\times\mathcal{X}$. The mutual information is in turn written as
\eqref{F-0.1} in Introduction as long as all the involved quantities (the
Shannon-Gibbs entropies of $X$, $Y$ and $(X,Y)$) are finite. This is nothing but
Voiculescu mentioned in \cite{V6} as an initial motivation of his introduction of
the liberation theory in free probability. Let us now apply the definition of
$I(X,Y)$ to our random matrix model of a given pair $(p,q)$ of projections. A random
matrix model at our disposal is the Grassmannian random matrix pair $(P(N),Q(N))$
whose joint distribution on $G(N,k(N))\times G(N,l(N))$ is the measure
$\lambda_N^{\psi_N}$ given in \eqref{F-2.3}. By the unitary invariance of trace
functions and the measure \eqref{F-1.2}, it is plain to see that the marginal
measures of $\lambda_N^{\psi_N}$ are $\gamma_{G(N,k(N))}$ and $\gamma_{G(N,l(N))}$;
thus  
$$
I(P(N),Q(N))  = S(\lambda_N^{\psi_N},\lambda_N^0)
\quad\mbox{with $\lambda_N^0=\gamma_{G(N,k(N))}\otimes\gamma_{G(N,l(N))}$}.
$$
In the proof of Theorem \ref{T-2.1} we obtained (see \eqref{F-2.8}) 
$$
-\chi_{\mathrm{proj}}(p,q) =
\lim_{N\rightarrow\infty}\frac{1}{N^2}
S(\lambda_N^{\psi_N},\lambda_N^0)
= \lim_{N\rightarrow\infty}\frac{1}{N^2} I(P(N),Q(N)). 
$$
Hence the minus free entropy of two projections can be also obtained as a scaling
limit of classical mutual information of our random matrix model. According to
\cite{V6} a ``heuristic definition" of {\it free mutual information} $i^*(p,q)$
should be
``$\chi_{\mathrm{proj}}(p) +\chi_{\mathrm{proj}}(q) - \chi_{\mathrm{proj}}(p,q)$"
on the analogy of classical theory. However, the actual definition of free mutual
information is completely different and based on the so-called liberation process
so that it may be particularly interesting to examine whether or not $i^*(p,q)$
coincides with $-\chi_{\mathrm{proj}}(p,q)$ in view of $\chi_\proj(p)=\chi_\proj(q)=0$.
In fact, our inequality in Theorem \ref{T-2.1} is a kind of logarithmic Sobolev
inequality and its right-hand side (the Dirichlet form part) is a ``derivative" of
$i^*(p,q)$, which also gives us a strong reason to do so. Let us give a heuristic
argument in the next subsection.

\subsection{$i^* = -\chi_{\mathrm{proj}}$ for two projections}
Although there are still some difficulties on regularity via the liberation
process, we give a heuristic computation indicating that $i^*(p,q)$ coincides with
$-\chi_\proj(p,q)$.

For a pair $(p,q)$ of projections in a tracial $W^*$-probability space
$(\cM,\tau)$, let $(\nu, \{\alpha_{ij}\}_{i,j=0}^1)$ be its representing data
consisting of a (not necessarily probability) measure $\nu$ on $(0,1)$ and the
trace values of four projections $E_{ij}$ (see \S\S1.1). Let $(u(t))_{t\geq0}$ be
a unitary free Brownian motion starting at $u(0)=\1$ (see \cite{Bi}). Letting
$p(t) := u(t)p u(t)^*$, a liberation process of projections starting
at $p$, we write $\nu_t$ and $E_{ij}(t)$ for $\nu$ and $E_{ij}$ corresponding to
$(p(t),q)$. Now, assume that $(\nu,\{\alpha_{ij}\}_{i,j=0}^1)$ satisfies the
assumptions of Proposition \ref{P-1.2}. By \cite[Corollary 8.6]{V6} the liberation
gradient $J_t := j(\bC p(t)+\bC(\1-p(t)):\bC q+\bC(\1-q))$ exists; hence
$\tau(E_{ij}(t)) =\alpha_{ij}$ for all $t \geq 0$ and $i,j=0,1$ thanks to
\cite[Lemma 12.5]{V6} together with $\tau(p(t)) = \tau(p)$. It is quite plausible
that each $\nu_t$ has the same properties as $\nu$, i.e., the assumptions of
Proposition \ref{P-1.2}. However, we could not derive these from the assumptions of
$\nu$ so that they have to be assumed here. Furthermore, we suppose that 
\begin{itemize} 
\item $f_t(x)$ is smooth in $(t,x) \in (0,+\infty)\times(0,1)$,
\end{itemize} 
which is also plausible to hold true. Set
\begin{gather*}
E(t) :=\1-(E_{00}(t)+E_{01}(t)+E_{10}(t)+E_{11}(t)), \\
X(t) := p(t)qp(t)+(1-p(t))q(1-p(t)), \\
\phi_t(x) := (Hf_t)(x) + \frac{\alpha_{01}+\alpha_{10}}{x}
-\frac{\alpha_{00}+\alpha_{11}}{1-x}\quad\mbox{for $0<x <1$},
\end{gather*}
where $Hf_t$ is the Hilbert transform of $f_t$. As stated in Proposition \ref{P-1.2},
the liberation gradient $J_t$ and the mutual free Fisher information
$\ffi^*(p(t) : q) :=
\ffi^*(\mathbb{C}p(t)+\mathbb{C}(1-p(t)) : \mathbb{C}q + \mathbb{C}(1-q))$
are given by    
\begin{gather*} 
J_t = [q,p(t)] \phi_t(E(t) X(t) E(t)), \\
\ffi^*(p(t):q) = \int_0^1 \phi_t(x)^2 x(1-x) f_t(x)\,dx.
\end{gather*} 
 
For each $m \in \mathbb{N}$ and $t,\varepsilon>0$ we set
$q(t):=u(\eps)u(t+\eps)^*qu(t+\eps)u(\eps)^*$ and
$J_t':=j(\bC p+\bC(\1-p):\bC q(t)+\bC(\1-q(t))$. Note \cite{Bi} that $u(\eps)$ is
$*$-free from $\{p,q(t)\}$ and that $(p,q(t))$ has the same distribution as
$(p,u(t)^*qu(t))$, which is clearly the same as $(p(t),q)$ too. Hence $(p,q(t),J_t')$
and $(p(t),q,J_t)$ behave in the same way under $\tau$. Therefore, by
\cite[Corollary 5.7]{V6} we have
\allowdisplaybreaks{
\begin{align*}
\tau\left(\left(p(t+\varepsilon)qp(t+\varepsilon)\right)^m\right)
&=
\tau\left(\left(p u(t+\varepsilon)^* q u(t+\varepsilon)\right)^m\right) \\
&=
\tau((u(\eps)pu(\eps)^*q(t))^m) \\
&=
\tau\left((pq(t))^m\right) + \frac{m\varepsilon}{2}
\tau\left([J'_t,p]\left(q(t) p q(t)\right)^{m-1}\right) + O(\varepsilon^2) \\
&=
\tau\left((p(t)qp(t))^m\right) + \frac{m\varepsilon}{2}
\tau\left([J_t,p(t)]\left(q p(t) q\right)^{m-1}\right) + O(\varepsilon^2)
\end{align*}
}so that
\begin{align*}
&{d\over d\eps}\bigg|_{\eps=0}
\tau\left(\left(p(t+\varepsilon)qp(t+\varepsilon)\right)^m\right) \\
&\qquad
={m\over2}\tau([J_t,p(t)](qp(t)q)^{m-1}) \\
&\qquad
={m\over2}\tau([[q,p(t)],p(t)](qp(t)q)^{m-1}\phi_t(E(t)X(t)E(t)) \\
&\qquad
=m\tau\left((p(t)qp(t))^m-(p(t)qp(t))^{m+1}\phi_t(E(t)(p(t)qp(t))E(t))\right)
\end{align*}
because both $E(t)$ and $X(t)$ are in the center of $\{p(t),q\}''$. In view of the
assumption on $(t,x) \mapsto f_t(x)$ the above equation implies that
\begin{align*}
{d\over dt}\int_0^1x^mf_t(x)\,dt
&=m\int_0^1(x^m-x^{m+1})\phi_t(x)f_t(x)\,dx \\
&=-\int_0^1x^m{\partial\over\partial x}(x(1-x)\phi_t(x)f_t(x))\,dx,
\end{align*}
which yields
\begin{equation}\label{F-3.1}
{\partial\over\partial t}f_t(x)
=-{\partial\over\partial x}(x(1-x)\phi_t(x)f_t(x)).
\end{equation}

Proposition \ref{P-1.1} says that
\begin{align*}
&\chi_\proj(p(t),q) \\
&\quad={1\over4}\int_0^1\int_0^1\log|x-y|\cdot f_t(x)f_t(y)\,dx\,dy \\
&\qquad+\frac{\alpha_{01}+\alpha_{10}}{2}\int_0^1\log x\cdot f_t(x)\,dx
+\frac{\alpha_{00}+\alpha_{11}}{2}\int_0^1\log(1-x)\cdot f_t(x)\,dx-C.
\end{align*}
Differentiating the above and applying \eqref{F-3.1} one can perform a heuristic
computation as follows:
\allowdisplaybreaks{
\begin{align*}
{d\over dt}\chi_\proj(p(t),q)
&=-{1\over2}\int_0^1\int_0^1\log|x-y|\cdot
{\partial\over\partial x}(x(1-x)\phi_t(x)f_t(x))f_t(y)\,dx\,dy \\
&\quad-{\alpha_{01}+\alpha_{10}\over2}\int_0^1\log x\cdot
{\partial\over\partial x}(x(1-x)\phi_t(x)f_t(x))\,dx \\
&\quad-{\alpha_{00}+\alpha_{11}\over2}\int_0^1\log(1-x)\cdot
{\partial\over\partial x}(x(1-x)\phi_t(x)f_t(x))\,dx \\
&={1\over2}\int_0^1\int_0^1{1\over x-y}\cdot x(1-x)\phi_t(x)f_t(x)f_t(y)
\,dx\,dy \\
&\quad+{\alpha_{01}+\alpha_{10}\over2}\int_0^1{1\over x}\cdot
x(1-x)\phi_t(x)f_t(x)\,dx \\
&\quad-{\alpha_{00}+\alpha_{11}\over2}\int_0^1{1\over1-x}\cdot
x(1-x)\phi_t(x)f_t(x)\,dx \\
&={1\over2}\int_0^1\phi_t(x)^2x(1-x)f_t(x)\,dx={1\over2}\ffi^*(p(t):q).
\end{align*}
}In particular, this shows that $\chi_{\mathrm{proj}}(p(t),q)$ is an increasing
function of $t\in(0,\infty)$. Moreover, since Theorem \ref{T-2.1} gives
$-\chi_{\mathrm{proj}}(p(t),q)\le\ffi^*(p(t):q)$ for all $t\ge0$, we have
$$
-\int_0^\infty\chi_\proj(p(t),q)\,dt
\le\frac{1}{2}\int_0^\infty\ffi^*(p(t):q)\,dt=i^*(p:q)<+\infty
$$
by \cite[Proposition 10.11.c]{V6} so that $\lim_{t\to\infty}\chi_\proj(p(t),q)=0$.
Therefore,
$$
i^*(p:q)={1\over2}\int_0^\infty\ffi^*(p(t):q)\,dt
=\Bigl[\chi_\proj(p(t),q)\Bigr]_0^\infty=-\chi_\proj(p,q)
$$
as long as $\chi_\proj(p,q)=\lim_{t\searrow0}\chi_\proj(p(t),q)$ is valid.

\subsection{A generalization of main theorem}
Free analogs of logarithmic Sobolev inequality were shown in \cite{Bi2} for the
single self-adjoint case and in \cite{HPU} for the single unitary case. The
inequality in Theorem 2.1 can be understood as such a free analog for two projections
with respect to the trivial ``hamiltonian" or trivial ``potential function." Thus it
is natural to make an attempt of generalizing it to the case where the given
``hamiltonian" is non-trivial. In fact, our method based on the random matrix
approximation still works for such an attempt. To do so, we need, on one hand, to give
the definitions of ``relative free entropy" and ``relative free Fisher information"
for pairs of projections, and on the other hand (as a technical side), to examine
Bakry and Emery's $\Gamma_2$-criterion by computing the Hessian of a certain trace
function on $G(N,k(N))\times G(N,l(N))$.

As mentioned in \cite[\S3]{HP1}, the distribution of a general pair of projections
can be understood as a tracial state on the $C^*$-algebra
\begin{equation*}
\mathcal{A} := \left\{ a(t) = [a_{ij}(t)] \in C([0,1]; M_2(\mathbb{C})) :
\text{$a(0), a(1)$ are diagonals} \right\}
\cong C^*(\mathbb{Z}_2\star\mathbb{Z}_2)
\end{equation*}
with the canonical generators of two projections: 
\begin{equation*} 
p(t) := \begin{bmatrix} 1 & 0 \\ 0 & 0 \end{bmatrix}\quad \text{and}\quad
q(t) := \begin{bmatrix} t & \sqrt{t(1-t)} \\ \sqrt{t(1-t)} & 1-t \end{bmatrix}.
\end{equation*}
The tracial state space of $\cA$ is denoted by $TS(\cA)$. An arbitrary
$\tau\in TS(\cA)$ is uniquely determined by the representing data
$(\nu,\{\alpha_{ij}\}_{i,j=0}^1)$ of $(p,q)$ in the GNS representation of $\cA$ with
respect to $\tau$ (see \S\S1.1); namely,
$$
\tau(a) = \alpha_{10}a_{11}(0) + \alpha_{01}a_{22}(0) + \alpha_{11}a_{11}(1)
+ \alpha_{00}a_{22}(1) + \int_0^1 \mathrm{tr}_2(a(t))\,d\nu(t)
$$
for every $a \in \mathcal{A}$, including \eqref{F-1.4}. We set
$\chi_\proj(\tau):=\chi_\proj(p,q)$ in the GNS representation with respect to $\tau$. 
Furthermore, let $h$ be a self-adjoint element in $\cA$ which we consider as a
general hamiltonian, and define a continuous function $\tilde h(t):=\Tr_2(h(t))$ on
$[0,1]$. Then, let us introduce the {\it relative free entropy}
$\widetilde{\Sigma}_h(\tau)$ of $\tau$ with respect to $h$ in the following way: When
$\alpha_{00}\alpha_{11}=\alpha_{01}\alpha_{10} = 0$, define
$$
\widetilde{\Sigma}_h(\tau):=-\chi_\proj(\tau)+\tau(h)+B_h,
$$
where $B_h$ is the maximum of the entropy functional
$\tau' \in TS(\mathcal{A}) \mapsto \chi_{\mathrm{proj}}(\tau') - \tau'(h)$ under the
condition that $\tau'(p) = \tau(p)$ and $\tau'(q)=\tau(q)$. More concretely,
\begin{align*} 
\widetilde{\Sigma}_h(\tau) &= -{1\over4}\Sigma(\nu)
+{1\over2}\int_0^1\bigl(\tilde{h}(x)
- (\alpha_{10}+\alpha_{01}) \log x \\
&\hskip3.5cm - (\alpha_{11}+\alpha_{00}) \log(1-x)\bigr)\,d\nu(x) + C_h,
\end{align*} 
where $C_h:=C+B_h$ with $C$ in \eqref{F-1.6}. Behind this definition is the large
deviation principle of the empirical distribution of the random projection matrix
pair $\lambda_N^h$ defined by replacing $\psi_N$ in \eqref{F-2.3} by $\tilde h$; its
proof is essentially same as in \cite{HP1}. Indeed, as (a) and (b) in the proof of
Theorem \ref{T-2.1}, we have
$C_h = \lim_{N\rightarrow\infty}\frac{1}{N^2}\log\widetilde{Z}_N^h$ for the
normalization constant $\widetilde{Z}_N^h$, and $\widetilde{\Sigma}_h(\tau)$ appears
as the rate function, justifying the term ``relative free entropy."

On the other hand, when $\tilde{h}(t)$ is assumed to be a $C^1$-function, let us
define the {\it relative free Fisher information} $\Phi_h(\tau)$ of $\tau$ with
respect to $h$ as follows: If the $\alpha_{ij}$'s satisfy the same condition as above
and moreover $\nu$ has the density $f := d\nu/dx \in L^3((0,1),x(1-x)dx)$, then
\begin{equation*}
\Phi_h(\tau) := \int_0^1 \bigg((Hf)(x)
+ \frac{\alpha_{10}+\alpha_{01}}{x} - \frac{\alpha_{11}+\alpha_{00}}{1-x}
- \tilde{h}'(x)\bigg)^2 x(1-x)\,d\nu(x);
\end{equation*}
otherwise $\Phi_h(\tau) := +\infty$. (Remark that the above integral is well-defined
permitting $+\infty$ and $\Phi_h(\tau)<+\infty$ is equivalent to the condition
\eqref{integrability}.) Note that when $\lambda_N^{\psi_N}$ is given under the same
assumptions (A) and (B) for $\tau$ as in the proof of Theorem \ref{T-2.1},
$\Phi_h(\tau)$ appears as the scaling limit of Dirichlet form:
\begin{equation*}
\lim_{N\rightarrow\infty} \frac{1}{2N^3} \int_{G(N,k(N))\times G(N,l(N))}
\left\Vert\nabla\log\frac{d\lambda_N^{\psi_N}}{d\lambda_N^h}(P,Q)
\right\Vert_{HS}^2\,d\lambda^{\psi_N}(P,Q).
\end{equation*}
In the trivial case where $h=0$, we have
$\widetilde\Sigma_h(\tau)=-\chi_\proj(p,q)$ and $\Phi_h(\tau)=\ffi^*(p:q)$ for $(p,q)$
in the GNS representation with respect to $\tau$ so that the next proposition is a
slight generalization of Theorem \ref{T-2.1}.

\begin{prop} Let $h$ be a self-adjoint element in $\mathcal{A}$. If $\tilde{h}(t)$
is a $C^2$-function on $[0,1]$ and it satisfies
$c_1\Vert\tilde{h}\Vert_{\infty}+c_2\Vert\tilde{h}''\Vert_{\infty} < 1$ for certain
universal constants $c_1, c_2 >0$, then the inequality 
\begin{equation*} 
\widetilde{\Sigma}_h(\tau) \leq \frac{1}{1-c_1\Vert\tilde{h}'\Vert_{\infty}
-c_2\Vert\tilde{h}''\Vert_{\infty}} \,\Phi_h(\tau)
\end{equation*} 
holds for every $\tau \in TS(\mathcal{A})$. 
\end{prop} 

\noindent
{\it Sketch of Proof.}\enspace
The proof is essentially same as that of Theorem 2.1, and the only difference is
in confirming Bakry and Emery's $\Gamma_2$-criterion \cite{BE} for $\lambda_N^h$
under the assumption of the proposition. The criterion in this case says that one
has a constant $c >0$ so that
$$
\mathrm{Ric}(G(N,k)\times G(N,l)) + \mathrm{Hess}(\Psi_N)
\geq c \cdot N I_{2k(N-k)+2l(N-l)}, 
$$ 
where $\mathrm{Hess}(\Psi_N)$ stands for the Hessian of the trace function
$$
\Psi_N(P,Q):=N\mathrm{Tr}_N(\tilde{h}(PQP))
\quad \mbox{for $(P,Q) \in G(N,k)\times G(N,l)$}.
$$
Thanks to \eqref{F-2.2} we need only to estimate $\mathrm{Hess}(\Psi_N)$ from below.
One can explicitly compute $\mathrm{Hess}(\Psi_N)$ in terms of the normal coordinate
\eqref{F-2.10}, which contains many terms of trace functions involving
$\tilde h'(PQP)$ and $\tilde h''(PQP)$. A rough estimation of the formula shows that
there are two universal constants $c_1, c_2 > 0$ so that 
\begin{equation*} 
\mathrm{Hess}(\Psi_N) \geq 
-N(c_1 \Vert\tilde{h}'\Vert_{\infty}+c_2 \Vert\tilde{h}''\Vert_{\infty})
I_{2k(N-k)+2l(N-l)},
\end{equation*}
while we do not know the best possible $c_1,c_2$. Now, the proposition follows from
the proof of Theorem \ref{T-2.1} together with the above estimate. \qed

\medskip
More details on computation of the Hessian $\mathrm{Hess}(\Psi_N)$ as well as the
constants $c_1,c_2$ in the above proof can be found in \cite[Remark 5.6]{HU}.

\end{document}